\documentclass[10pt,twocolumn]{IEEEtran}
 \usepackage{graphicx}
 \usepackage{amsmath}
 \usepackage{amsfonts}
\newtheorem{defin}{Definition}
\newtheorem{theorem}{Theorem}
\newtheorem{lemma}{Lemma}
\newtheorem{problem}{Problem}
\newtheorem{remark}{Remark}
\newenvironment{prf}{\smallbreak\noindent{\it Proof: }}{\hfill$\Box$\smallbreak}

\DeclareMathOperator{\dist}{dist}
\DeclareMathOperator{\sos}{SoS}

\begin{document}
\title{Transverse Dynamics and\\ Regions of Stability for\\ Nonlinear Hybrid Limit Cycles 
\thanks{This work was supported by NSF Contract 0915148.}}
\author{Ian R. Manchester\\ 
CSAIL, Massachusetts Institute of Technology, irm@mit.edu}       
\maketitle

\begin{abstract}
This paper presents an algorithm for computing inner estimates of the regions of attraction of limit cycles of a nonlinear hybrid system. The basic procedure is: (1) compute the dynamics of the system transverse to the limit cycle; (2) from the linearization of the transverse dynamics construct a quadratic candidate Lyapunov function; (3) search for a new Lyapunov function verifying maximal regions of orbital stability via iterated of sum-of-squares programs. The construction of the transverse dynamics is novel, and valid for a broad class of nonlinear hybrid systems. The problem of stabilization of unstable limit cycles will also be addressed, and a solution given based on stabilization of the transverse linearization.
\end{abstract}

\section{Introduction}

Nonlinear dynamical systems exhibiting oscillating solutions are found in an extraordinary variety of engineering and scientific problems. For example: cellular signalling, radio-frequency circuits, walking robots, and population dynamics. Stability analysis of such oscillations has a long history, with the modern theory going back to Poincar\'e. For planar systems substantial qualitative insight can be gained, however for higher-order systems almost all stability results are local. In this paper, we aim to characterize {\em regions of stability} to limit cycles of nonlinear systems. The problem of {\em orbital stabilization} via feedback control will also be addressed.

A major motivation for the work in this paper is control of underactuated ``dynamic walking'' robots (\cite{Collins05}). These robots can exhibit efficient, naturalistic, and highly dynamic gaits. However, control design and stability analysis for such robots is a challenging task since their dynamics are intrinsically hybrid and highly nonlinear. Local stabilizing control design has been investigated via hybrid zero dynamics (\cite{Westervelt03a, Westervelt07b}) and transverse linearization (see, e.g., \cite{Shiriaev08, Shiriaev09, Manchester09}).

As well as being of interest in their own right, estimates of regions of stability will be an enabling technology for planning transitions among a library of stabilized walking gaits (\cite{Gregg10}), and for constructive control design and motion-planning (\cite{Tedrake10}).

The most well-known tool for analysis of limit cycles is the Poincar\'e map: orbital stability is characterized by stability of an associated ``first-return map'', describing the repeated passes of the system through a single transversal hypersurface. Often a linearization of the first-return map is computed numerically, and its eigenvalues can be used to verify local orbital stability. Since the system's evolution is only analyzed on a single surface, regions of stability in the full state-space are difficult to evaluate via the Poincar\'e map.

A related technique known variably as ``transverse coordinates'' or ``moving Poincar\'e sections'' also has a long history and was certainly known to exist by Poincar\'e, however has not been much used in applications until recently due to difficulty in the relevant computations (see \cite{Hale80}). With this technique, a new coordinate system is defined on a family of transversal hypersurfaces which move about the orbit under study. In most cases, it is also used to study local stability, however as we will show it can be adapted to characterize regions of stability in the full state space.

The construction of the transverse dynamics is novel, analytical, valid for a large class of systems, and amenable to analysis via sum-of-squares (SoS) programming. Previous work by the author and colleagues has utilized a construction specifically for Lagrangian mechanical systems (see \cite{Shiriaev08, Shiriaev09, Manchester09}). The new construction is also useful for design of stabilizing controllers for non-periodic trajectories of highly nonlinear systems without this special structure, e.g. in \cite{Shkolnik10a} a quadruped robot bounding over rough terrain was stabilized using a preliminary version of the construction in the present paper.

In this paper we present an algorithm to compute a conservative estimate of the region of stability to a known periodic solution of hybrid nonlinear system. The method we propose is to construct the transverse dynamics in regions of the orbit, and then  utilize the well-known SoS relaxation of polynomial positivity which is amenable to efficient computation via semidefinite programming  (see, e.g., \cite{Shor87, Parrilo00, Parrilo03a}). The sum-of-squares relaxation has been previously used to characterize regions of stability of equilibria of nonlinear systems (see, e.g., 
\cite{Topcu08, Tan08, Henrion05}). 

The method we propose has been tested on a number of example systems, presented in the companion paper \cite{Manchester10b}. There is comparatively little work on computing regions of stability of limit cycles. The proposed method has aspects in common with the surface Lyapunov functions proposed in \cite{Goncalves05}, however that method was restricted to piecewise linear systems.  The technique of cell-to-cell mapping, proposed by  \cite{Hsu80},  improves the efficiency of exhaustive grid-based methods of regional analysis and has been used in analysis of walking robots (\cite{Schwab01}), however the computational cost is still exponential in the dimension of the system.


\section{Preliminaries and Problem Statements}

\subsection{Stability of Limit Cycles}
Consider an autonomous dynamical system, which may be continuous or hybrid, with a state $x\in \mathbb R^n$. The solution, or {\em flow} of the system is denoted by
$\Phi(x_0,t)$, i.e. $x(t) = \Phi(x_0, t)$ is the solution at time $t>0$ of the dynamical system in question from an initial state $x(0) = x_0$. If this solution exists and is unique then $\Phi(x_0, t)$ is well-defined. Note that since the system is autonomous, we can consider flows starting from $t=0$ without any loss of generality.

Suppose the system has a non-trivial $T$-periodic orbit, i.e. $T>0$ is the minimal period such that $x^\star(t)=x^\star(t+T)$ for all $t$, and one would like to analyze the stability of this orbit. It is well-known that such a solution cannot be asymptotically stable in the standard sense, since perturbations in phase are persistent. The more appropriate notion is {\em orbital stability}. The definitions in this section are all standard (see, e.g., \cite{Hale80, Hauser94}).

\begin{defin} Consider non-trivial $T$-periodic solution $x^\star(t)$ of a dynamical system with flow $\Phi(\cdot, \cdot)$, and let $\Gamma^\star$ denote the solution curve: $\Gamma^\star = \{x\in\mathbb R^n : \exists t\in [0,T) : x = x^\star(t)\}$. The solution $x^\star(\cdot)$ is said to be {\em asymptotically orbitally stable} if there exists a $b>0$ such that for any $x_0$ satisfying $\dist(x_0, \Gamma^\star)<b$ the solution exists, is unique, and $\dist(\Phi(x_0,t),\Gamma^\star)\rightarrow 0$ as $t\rightarrow \infty$.
\end{defin}
The distance to a set is defined in the usual way: $\dist(x, \Gamma^\star) = \inf_{y\in \Gamma^\star}{|y-x|}$ with $|\cdot|$ the Euclidean norm in $\mathbb R^n$.

A stronger statement is {\em exponential orbital stability}:
\begin{defin} A $T$-periodic solution $x^\star(\cdot)$ is said to be {\em exponentially orbitally stable} if it is orbitally stable and furthermore there exists a $b>0, K>0, c>0$ such that for any $x_0$ satisfying $\dist(x_0, \Gamma^\star)<b$ we have 
\[
\dist(\Phi(x_0,t),\Gamma^\star)\le K\dist(x_0,\Gamma^\star)e^{-ct}.
\]
\end{defin}

The primary aim of this paper is to characterize regions of orbital stability:
\begin{defin} A set $R\subset \mathbb R^n$ with non-empty interior $R^i$ and  $\Gamma^\star\subset R^i$ is said to be an {\em inner estimate of the region of stability of }$x^\star(\cdot)$ if for all $x_0\in R$ we have $\dist(\Phi(x_0,t),\Gamma^\star)\rightarrow 0$ as $t\rightarrow \infty$.
\end{defin}
\subsection{Problem Statements}

This paper will suggest algorithms for three analysis and control problems for limit cycles:
\begin{problem}\label{prob:regions}
Given an autonomous smooth nonlinear system with state $x\in \mathbb R^n$:
\begin{equation}\label{eqn:sys}
\dot x = f(x)
\end{equation}
and a non-trivial $T$-periodic solution $x^\star(t)=x^\star(t+T)$ of \eqref{eqn:sys} such that $f(x^\star(t))\ne 0$ for any $t \in [0,T]$. Characterize the stability of $x^\star(\cdot)$ and if it is exponentially orbitally stable then compute an inner estimate of its region of stability. \hfill$\Box$\smallbreak
\end{problem}

Many systems are best modelled by a combination of smooth nonlinear dynamics with occasional moments of instantaneous change in the state. This may model a discrete change in a switching controller, or a period of extremely fast change in the system state, e.g. an impact event in a physical system.

The second problem is to estimate regions of stability for a class of hybrid nonlinear systems:
\begin{problem}\label{prob:regions_hybrid}
Consider an autonomous hybrid nonlinear system with state $x\in \mathbb R^n$ and switching dynamics  defined between hyperplanes:
\begin{eqnarray}
\dot x &=& f(x),  \ x\notin S^-\label{eqn:sysh1}\\
x^+&=&\Delta(x), \ x\in S^-. \label{eqn:sysh2}
\end{eqnarray}
Suppose $f(\cdot)$ and $\Delta(\cdot)$ are smooth and $\Delta:S^-\rightarrow S^+$ where
\begin{eqnarray}
S^- &=& \{x: c_-'x=d_-, g(x)\ge0\},\\
S^+ &=& \{x: c_+'x=d_+\},
\end{eqnarray}
$c_-, c_+\in \mathbb R^n$, and $d_-, d_+ \in \mathbb R$. Suppose $x^\star(\cdot)$ is a non-trivial $T$-periodic solution that undergoes $N$ impacts at times $\{t_1, t_2, ..., t_N\}+kT$ for integer $k$. 
We will assume that the impacts are not ``grazing'', i.e. $c_-'f(x^\star(t_i)) \ne 0$ and $c_+'f(x^\star(t_i)) \ne 0 $ for all $i$.

The problem statement is to characterize the stability of $x^\star(\cdot)$ and if it is exponentially orbitally stable compute an inner estimate of its region of stability. \hfill$\Box$\smallbreak
\end{problem}
For simplicity of expression we will consider the problem with a single switching surface and a single set of continuous dynamics, however the extension to multiple switches and multiple continuous phases is trivial. It is also straightforward to have the dimension of the continuous system change between different phases. 

The main practical restriction in this class is that switching surfaces are planar. This is quite a strong restriction, but it greatly simplifies proving orbital stability via planar transversal surfaces, since we can make the transversal surface line up with the switching surfaces before and after impact. Furthermore, it is true for some important and common models of walking robots such as the rimless wheel and the compass gait. For other systems it may be possible to contstruct a change of coordinates such that the impact map is planar, e.g. for a multi-link walking robot one can choose one generalized coordinate to be the distance from the swing foot to the ground plane.

The third problem is one of feedback orbital stabilization:
\begin{problem}\label{prob:stabilization}
Consider a controlled, possibly hybrid, system with a control input $u\in \mathbb R^m$:
\begin{eqnarray}\label{eqn:sys_u1}
\dot x &=& f(x,u),  \ x\notin S^-\\
x^+&=&\Delta(x), \ x\in S^-.
\end{eqnarray}
Suppose this system has a $T$-periodic solution $x^\star(\cdot)$ that undergoes $N$ impacts per period at times $\{t_1, t_2, ..., t_N\}+kT$ for integer $k$, and is generated by a piecewise continuous control signal, $u^\star(t) = u^\star(t+T)$. If this solution is not exponentially orbitally stable then, if possible, construct a state-feedback exponentially orbitally stabilizing controller and compute an inner estimate of its region of orbital stability. \hfill$\Box$\smallbreak
\end{problem}

\section{Regions of Stability for Continuous Systems}
In this section we propose a solution to Problem \ref{prob:regions}.

The process we propose for finding regions of orbital stability is based on the construction of a smooth local change of coordinates $x \rightarrow (x_\perp, \tau)$. At each point $t\in [0, T]$ we define a hyperplane $S(t)$, with $S(0)=S(T)$, which is transversal to the solution $\Gamma^\star$, i.e. $\dot x^\star(t)\not\in S(t)$. 

\begin{figure}[t]
\centering
\hskip -10pt\includegraphics[width=0.65\columnwidth]{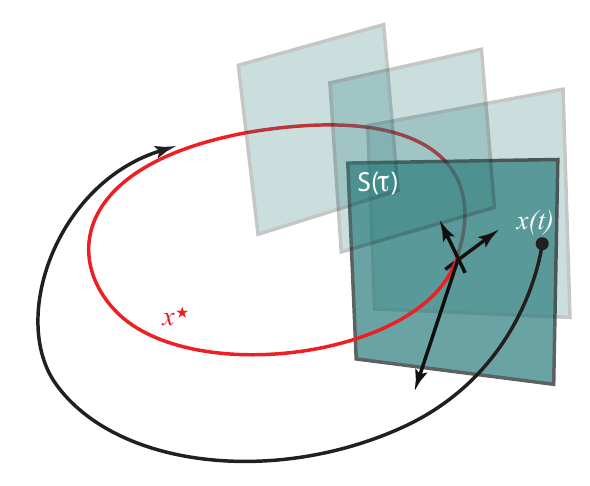}\\
\hskip 20pt\includegraphics[width=0.66\columnwidth]{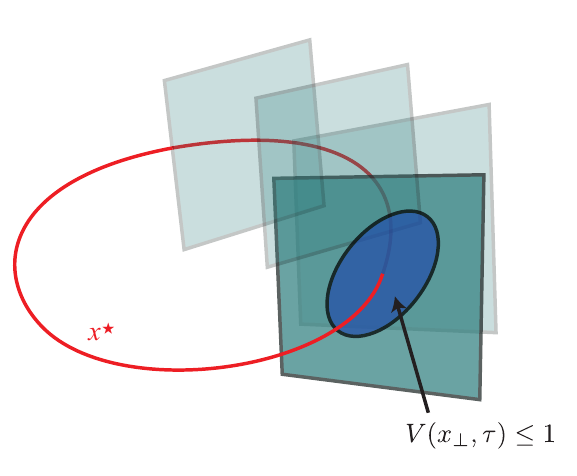}
\caption{Left: transversal surfaces $S(\tau)$ around the target orbit $x^\star$, with a particular solution $x(t)$ converging to the orbit $x^\star$. Right: a Lyapunov function defined on a transversal surface.}
\label{fig:ts1}
\end{figure}

Given a point $x$ nearby $x^\star(\cdot)$, the scalar $\tau \in [0, T)$ represents which of these transversal sufaces $S(\tau)$ the current state $x$ inhabits; the vector $x_\perp \in \mathbb R^{n-1}$ is the ``transversal'' state representing the location of $x$ within the hyperplane $S(\tau)$, with $x_\perp = 0$ implying that $x=x^\star(\tau)$. This is visualised in Figure \ref{fig:ts1}. 

We will show that in some region around the target orbit $\Gamma^\star$, the dynamics in the new coordinate system are well-defined and have the form:
\begin{eqnarray}
\dot x_\perp &=&A(\tau)x_\perp+h(x_\perp, \tau),\label{eqn:trans_full1}\\
\dot \tau &=& 1+g( x_\perp, \tau), \label{eqn:trans_full2}
\end{eqnarray}
with $h=O(|x_\perp|^2)$ and $g=O(|x_\perp|)$ near $x_\perp=0$ for all $\tau$. The fact that there is a differential equation for $\tau$ corresponds to the following fact: if the state $x(t_1) \in S(\tau_1)$, then after a short time interval $\delta_t$ it does {\em not} necessarily follow that $x(t_1+\delta_t) \in S(\tau_1+\delta_t)$.

The main reason for this decomposition is that if we can prove that $x_\perp \rightarrow 0$ then we have proven orbital stability of $x^\star$.

Associated with the above nonlinear system is the {\em transverse linearization}: a first-order approximation of the $x_\perp$ dynamics expressed as a periodic linear system:
\begin{equation}\label{eqn:translin}
\dot x_\perp =A(t)x_\perp
\end{equation}
for $t\in [0,T)$.

It is known that the periodic solution $x^\star$ of the nonlinear system $\dot x = f(x)$ is exponentially orbitally stable if and only if the periodic linear system \eqref{eqn:translin} is asymptotically and hence exponentially stable (see, e.g., \cite{Hale80, Hauser94}).

The process for computing regions of orbital stability is as follows:
\begin{enumerate}
\item Select a family of transversal surfaces $S(\tau)$, and the associated transformation $x\rightarrow (x_\perp, \tau)$.
\item Compute the nonlinear dynamics in this new coordinate system as well as a periodic linear system representing the dynamics of $x_\perp$ close to the orbit: the {\em transverse linearization}.
\item Construct a candidate quadratic Lyapunov function associated with the transverse linearization via standard techniques from linear control theory.
\item Using this result as an initial seed, iteratively solve a sequence of sum-of-squares programs to compute maximal regions in which a Lyapunov function can be found verifying both well-posedness of the change of coordinates and orbital stability for the true nonlinear dynamics.
\end{enumerate}
The details of each step are given in the following four subsections

\subsection{Selection of a Set of Transversal Surfaces}

Suppose we have a periodic orbit $x^\star(t) = x^\star(t+T)$ with $\dot x^\star(t)\ne 0 \, \forall t\in[0,T)$. At each point $x^\star(\tau)$ of the target orbit we define a transversal surface $S(\tau)$ in the following way:
\[
S(\tau)=\{y\in\mathbb R^n : z(\tau)'(y-x^\star(\tau))=0\}
\]
where $z(\tau):[0, T)\rightarrow \mathbb R^n$ is a smooth periodic vector function to be chosen. We will enforce that $z(t)$ has bounded derivative. In the literature on the use of transversal coordinates to prove local properties about periodic solutions, it is common to choose $z(\tau)=f(x^\star(\tau))$ \cite{Hale80}, \cite{Hahn67}. That is, the transversal planes are orthogonal to the current motion of the system.

However, orthogonal transversal planes are often a bad choice when performing analysis on larger regions around the orbit, and allowing some freedom in $z(\tau)$ is highly beneficial. The reason is, roughly speaking, that singularities will occur in the change of coordinates $x\rightarrow (x_\perp, \tau)$ near sections of $x^\star$ with large curvature. This will be made more precise in Section \ref{sec:choose_z}.

The primary requirement is that the resulting $S(\tau)$ are still {\em transversal} to the orbit, which is guaranteed if there is some $\delta>0$ such that $z(\tau)'f(x^\star(\tau))>\delta$ for all $\tau\in [0, T)$. I.e., $z(\tau)$ is never orthogonal to $f(x^\star(\tau))$.

As a technical condition, we require that $z(\tau)$ be Lipschitz on each continuous interval. For simplicity of derivations, and without loss of generality, we will further enforce that $\|z(\tau)\|=1$ for all $\tau$. If planes orthogonal to the system motion are desired, we can take $z(\tau) = f(x^\star(\tau))/\|f(x^\star(\tau))\|$.

\subsection{Construction of a Moving Coordinate System}

Having chosen a set of hyperplanes, defined by $z(\tau)$, we now construct a smoothly $\tau$-varying coordinate system upon this subspace.
For $n=2$, construction of a basis is trivial: pick, e.g., $[-z_2, \ z_1]'$ as the basis vector for $\Pi$. For $n\ge 3$, a constructive can be adapted from a method in \cite{Diliberto56} and \cite[Ch. VI]{Hale80}.

\begin{enumerate}
\item Choose a vector $w$ such that $w$ is not collinear with  $z(\tau)$ for any $\tau$. 
\item Choose a fixed orthonormal basis $\eta_j, j = 1, 2, ..., n$ for $\mathbb R^n$ with $w$ as its first element.
\item for each $\tau\in[0,T)$, define the plane containing both $w$ and $z(\tau)$, and the rotation matrix $R(\tau)$ which takes $w$ to be collinear with $z(\tau)$ rotating in this plane.
\item Define $\xi_j(\tau) = R(\tau)\eta_j$ then $\xi_2(\tau), ... , \xi_n(\tau)$ form a basis for the transversal coordinates.
\end{enumerate}

An explicit formula for $\xi_j(\tau)$ in terms of $\eta_1$, $\eta_j$ and $z(\tau)$ is
\begin{equation}\label{eqn:xi}
\xi_j(\tau) = \eta_j-\frac{\eta_j'z(\tau)}{1+\eta_1'z(\tau)}(\eta_1+z(\tau)), \ j = 2, 3, ..., n.
\end{equation}
Note that since $\eta_1$ and $z(\tau)$ are, by construction, unit vectors which are not collinear, $\eta_1'z(\tau)<1$ for all $\tau$ so \eqref{eqn:xi} is well-defined.

The following lemma suggests that in practice the vector $w$ in the above construction can be chosen at random:
\begin{lemma}\label{lem:w}
Suppose $z(t):[0,T]\rightarrow \mathbb{R}^n, n\ge 3$ is piecewise continuous with a finite number of points of discontinuity $t_i$, and Lipschitz on each interval $[t_i, t_{i+1})$. Then the set of unit vectors $w$ such that $w=\pm z(t)$ for some $t\in[0,T]$ has measure zero.
\end{lemma}
\begin{prf}
For each interval $[t_i, t_{i+1})$ consider a closed interval $[t_i, t_{i+1}]$ with $z(t_{i+1})$ temporarily defined as 
\[
z(t_{i+1}):=\lim_{t<t_{i+1}, t\rightarrow t_{i+1} } z(t)
\]
Now, since $z(t)$ is unit-length, we can consider the set $z(t), t \in [t_i, t_{i+1}]$ to be a curve $\mathcal Z_i$ on the unit sphere $\mathbb R^n$. Since $z(t)$ has bounded first derivative, it is Lipschitz on $[t_i, t_{i+1}]$, and hence the curve $\mathcal Z_i$ is rectifiable, and it is known that a rectifiable curve on a unit sphere in $\mathbb R^n, n\ge 3$ covers a set of finite measure.

Since there is a finite number of continuous intervals, the union $\bigcup_i \mathcal Z_i$ is clearly also of zero measure.
\end{prf}

Having this basis defined, we also construct the projection operator:
\[
\Pi(\tau) = \begin{bmatrix}\xi_2(\tau)'\\ \vdots \\ \xi_n(\tau)'\end{bmatrix}
\]
which defines the mapping $x\rightarrow (x_\perp, \tau)$. That is, if $x\in S(\tau)$ then 
\[
x_\perp = \Pi(\tau)x.
\]
Note that a given $x\in \mathbb R^n$ will in general be in more than one transversal plane, i.e. we can have $x\in S(\tau_1)$ and $x\in S(\tau_2)$ with $\tau_1 \ne \tau_2$, though this will not cause any problems for the proposed method.


\subsection{Transverse Dynamics and Linearization}

\begin{theorem}
The dynamics of the system in the new coordinates $(x_\perp, \tau)$ are given by:
\begin{eqnarray}
\dot x_\perp&=& \dot \tau\left[\frac{d}{d\tau}\Pi(\tau)\right]\Pi(\tau)'x_\perp+\Pi(\tau)f(x^\star(\tau)+\Pi(\tau)'x_\perp)\notag\\ &&\,\,\,-\Pi(\tau)f(x^\star(\tau)) \dot \tau, \label{eqn:xperpdot}\\
\label{eqn:taudot}
\dot \tau &=& \frac{z(\tau)'f(x^\star(\tau)+\Pi(\tau)'x_\perp)}{z(\tau)'f(x^\star(\tau))-\frac{d z(\tau)}{d \tau}'\Pi(\tau)'x_\perp}.
\end{eqnarray}
\end{theorem}

\begin{prf}
Consider the transversal coordinate $x_\perp(t) = \Pi(\tau)(y(t)-x^\star(\tau))$ where $\tau$ is such that $z(\tau)'(y-x^\star(\tau))=0$.

\[
\dot x_\perp = \dot \tau\left[\frac{d}{d\tau}\Pi(\tau)\right][y(t)-x^\star(\tau)]+\Pi(\tau)\left[f(y)-f(x^\star(\tau)) \dot \tau\right]
\]
Since $y= x^\star(\tau)+\Pi(\tau)'x_\perp$ this can be written as \eqref{eqn:xperpdot}

To find the dynamics of $\tau$, consider the orthogonality condition:
\begin{equation}\label{eqn:orthog}
F(t,\tau) = z(\tau)'(x(t)-x^\star(\tau))=0.
\end{equation}
Since this remains true under evolution of the system, the total derivative of $F(t,\tau)$ is zero:
$
\dot F(t, \tau)= \frac{\partial F}{\partial \tau}\dot \tau +\frac{\partial F}{\partial t}=0,
$
and, by the implicit function theorem, if $\frac{\partial F}{\partial \tau}\ne 0$ we have
$
\dot \tau = -\left(\frac{\partial F}{\partial \tau}\right)^{-1}\frac{\partial F}{\partial t},
$
after some straightforward manipulations from \eqref{eqn:orthog}, we find that
\begin{equation}\label{eqn:taudoty}
\dot \tau = \frac{z(\tau)'f(x(t))}{z(\tau)'f(x^\star(\tau))-\frac{\partial z(\tau)}{\partial \tau}'(x(t)-x^\star(\tau))}.
\end{equation}
As above, this can be given in terms of $\tau$ and $x_\perp= \Pi(\tau)(x(t)-x^\star(\tau))$:
\begin{equation}\label{eqn:taudot}
\dot \tau = \frac{z(\tau)'f(x^\star(\tau)+\Pi(\tau)'x_\perp)}{z(\tau)'f(x^\star(\tau))-\frac{\partial z(\tau)}{\partial \tau}'\Pi(\tau)'x_\perp}.
\end{equation}

\end{prf}

\subsubsection{Transverse Linearization}

The transverse linearization is the following $(n-1)$-dimensional $T$-periodic linear system:
\begin{equation}\label{eqn:translin1}
\dot z = A(t) z
\end{equation}
where $A(t)$ comes from differentiating \eqref{eqn:xperpdot} with respect to $x_\perp$:
\begin{eqnarray}
A(t) &=& \left[\frac{d}{dt}\Pi(t)\right]\Pi(t)'+\Pi(t) \frac{\partial f(x^\star(t))}{\partial x}\Pi(t)'\notag\\ &&-\Pi(t)f(x^\star(t)) \left.\frac{\partial\dot\tau}{\partial x_\perp}\right|_{x_\perp=0}\label{eqn:translinA}
\end{eqnarray}
where
\begin{equation}\label{eqn:dtaudxp}
\left.\frac{\partial\dot\tau}{\partial x_\perp}\right|_{x_\perp=0}= \frac{z(t)'\frac{\partial f(x^\star(t))}{\partial x}\Pi(t)' + \frac{\partial z(t)}{\partial t}'\Pi(t)'}{z(t)'f(x^\star(t))}.
\end{equation}

\begin{remark} In two special cases some simplifications are possible: if the system is planar, i.e. $n=2$, then $[\frac{d}{d\tau}\Pi(\tau)]\Pi(\tau)'=0$ so \eqref{eqn:xperpdot} and \eqref{eqn:translinA} can be simplified by removing the first term from each. If transversal planes orthogonal to the system motion are chosen then $\Pi(\tau)f(x^\star(\tau)) =0$ so \eqref{eqn:xperpdot} and \eqref{eqn:translinA} can be simplified by removing the last term from each.
\end{remark}

\begin{remark}
Note that in the selection of $z(\tau)$ we have enforced that $z(\tau)'f(x^\star(\tau))>\delta$ for all $\tau$, so if $x_\perp$ remains sufficiently small we will have:
\[
z(\tau)'f(x^\star(\tau))-\frac{d z(\tau)}{d \tau}'\Pi(\tau)'x_\perp\ge \alpha >0
\]
for some $\alpha$, so the dynamics of $\tau$ are well-defined. This condition breaks if there is a continuum of $\tau$ such that $x^\star(\tau)$ satisfies $f(x^\star(\tau))'y-x^\star(\tau)=0$. E.g. for constant speed curves ($|\dot x|=0$) this means $|y-x^\star(\tau)|$ is exactly the radius of curvature of the the target orbit, i.e.  $1/ |\ddot x|$. In Section \ref{sec:choose_z} we will discuss optimizing $z(\tau)$ so as to maximize the regions in which the dynamics of $\tau$ are well-defined.
\end{remark}

\subsection{Verification of Orbital Stability}

We now state conditions which guarantee regions of stability, giving a solution to Problem \ref{prob:regions}. It will then be shown how to optimize regions for polynomial systems via a SoS relaxation.

\begin{theorem}\label{thm:stab}
Suppose there exists a function $V:\mathbb R^{n-1}\times \mathbb R \rightarrow \mathbb R$ such that $D:=\{(x_\perp, \tau):V(x_\perp, \tau)\le 1\}$ is compact and for which following inequalities hold for all $(x_\perp, \tau)\in D$:
\begin{eqnarray}
V(x_\perp,\tau)&>&0, \  \  x_\perp \ne 0\label{eqn:v}\\
\frac{d}{dt}V(x_\perp, \tau)&<&0, \  \  x_\perp \ne 0\label{eqn:vdot}\\
V(0, \tau)=\frac{d}{dt}V(0, \tau)&=&0,\label{eqn:v:vd:0}\\
z(\tau)'f(x^\star(\tau))-\frac{\partial z(\tau)}{\partial \tau}'\Pi(\tau)'x_\perp&>&0.\label{eqn:well-posed}
\end{eqnarray}
then the $D$ is an inner estimate for the region of orbital stability of $x^\star(\cdot)$.
\end{theorem}

\begin{prf}  It follows from the smoothness assumptions on $f(x)$ and \eqref{eqn:well-posed} that both $\dot x_\perp$ and $\dot \tau$ are Lipschitz in $D$, and by the usual Lyapunov arguments \eqref{eqn:v}, \eqref{eqn:vdot} ensure that the compact set $D$ is invariant. Therefore the solutions of $\tau$ and $x_\perp$ from any initial conditions in $D$ exist and are unique. Furthermore, \eqref{eqn:vdot}, \eqref{eqn:v:vd:0} imply that from any initial condition in $D$ we have $x_\perp \rightarrow 0$, from which it follows that $x^\star$ is orbitally asymptotically stable and $D$ is an inner estimate of its region of stability.
\end{prf}

The optimization problem is then to search over Lyapunov functions $V(x_\perp, \tau)$ so as to maximize the size of the regions $D$ such that these conditions can be verified.

In practice, we will often strengthen the above statements by choosing some small constants $\delta_1>0, \delta_2>0, \delta_3>0$ and requiring:
\begin{eqnarray}
V(x_\perp, \tau) -\delta_1|x|^2&\ge&0, \label{eqn:strict_Vdot}\\
\frac{d}{dt}V(x_\perp, \tau) +\delta_2|x|^2&\le&0, \label{eqn:strict_Vdot}\\
z(\tau)'f(x^\star(\tau))-\frac{\partial z(\tau)}{\partial \tau}'\Pi(\tau)'x_\perp-\delta_3&\ge&0.  \label{eqn:strict_well_posed}
\end{eqnarray}
This improves the robustness of the numerical procedures for verification.


\subsubsection{Sums-of-Squares Relaxation}

\ 

We now show how the regions can be computed using a SoS relaxation. For the purposes of this section, we assume that the vector field $f(x)$ is polynomial in $x$. This being the case, let us fix a particular value of $\tau$ and examine the formula for $\dot \tau$:
\[
\dot \tau =\frac{z(\tau)'f(x^\star(\tau)+\Pi(\tau)'x_\perp)}{z(\tau)'f(x^\star(\tau))-\frac{\partial z(\tau)}{\partial \tau}'\Pi(\tau)'x_\perp}=:\frac{n(x_\perp, \tau)}{d(x_\perp, \tau)}.
\]
It is clear that both the numerator and denominator,  $n(x_\perp, \tau)$ and $d(x_\perp, \tau)$ respectively, are both polynomial in 
$x_\perp$. 

Furthermore, the well-posedness condition ensures that $d(x_\perp,\tau)>0$, hence we can multiply both sides of condition \eqref{eqn:strict_Vdot} by $d(x_\perp, \tau)$ without changing its validity, resulting in the condition $DV(x_\perp, \tau)\le 0$, where $DV(x_\perp,\tau)$ is given by \eqref{eqn:stab1}. It is straightforward to verify that this condition is also polynomial in $x_\perp$.

\begin{figure*}[!t]
\begin{eqnarray}
DV(x_\perp, \tau)&:=&\frac{\partial}{\partial \tau}V(x_\perp, \tau)n(x_\perp, \tau)\notag + \frac{\partial}{\partial x_\perp}V(x_\perp, \tau) \notag\Bigg [ n(x_\perp, \tau)  \tau\left[\frac{d}{d\tau}\Pi(\tau)\right]\Pi(\tau)'x_\perp\notag\\ &&+d(x_\perp, \tau) \Pi(\tau)f(x^\star(\tau)+\Pi(\tau)'x_\perp)  -\Pi(\tau)f(x^\star(\tau)) n(x_\perp, \tau)\Bigg]+d(x_\perp, \tau) \delta_1|x_\perp|^2. \label{eqn:stab1}
\end{eqnarray}
\hrulefill
\vspace*{4pt}
\end{figure*}

We can verify these conditions regionally using Lagrange multipliers $l(x_\perp)$ and $m(x_\perp)$  that are polynomial in $x_\perp$ with the following SoS constraints:
\begin{eqnarray}
-DV(x_\perp,\tau)-l(x_\perp)(1-V(x_\perp, \tau)) &=& \sos,\label{eqn:sos1}\\
d( \tau, x_\perp)-\delta_2-m(x_\perp)(1-V(x_\perp, \tau))&=& \sos,\label{eqn:sos2}\\
l(x_\perp)  &=& \sos,\label{eqn:sos3}\\
m(x_\perp)  &=& \sos. \label{eqn:sos4}
\end{eqnarray}

In practice we sample a sufficiently fine finite sequence $\tau_i, i=1, 2, ..., N_i$ such that $\tau_1 = 0, \tau_{N_i}=T$. Then for each $\tau_i$ and verify the conditions \eqref{eqn:sos1} -- \eqref{eqn:sos4} for each fixed $\tau_i$.

The objective is to maximize the regions satisfying \eqref{eqn:sos1} -- \eqref{eqn:sos4}. The decision parameters are the coefficients of $V$ (within some restricted class) and the Lagrange multipliers $l$ and $m$ at each sample of $\tau$. This optimization is bilinear in the decision parameters, however a reasonable approach which has proven successful is to iterate between optimizing over multipliers and optimizing over $V$.

If $V$ is a polynomial of higher order than quadratic, then what exactly should be optimized in the search for $V$ is open to some choice, but a natural candidate is to maximize the size of some ball in $\|x_\perp\|$ contained in the set $V(x_\perp, \tau)\le 1$. This can be expressed as a sum-of-squares program, and has been found to be effective for computing regions of stability for equilibria (see \cite{Topcu08}).

\subsection{An Initial Candidate Lyapunov Function}

The iteration approach to solving a bilinear optimization problem requires a reasonably good initial guess. We propose using the solution of the periodic Lyapunov differential equation for the transverse linearization as an initial guess for $V(x_\perp, \tau)$.

It is well known that the transverse linearization, a periodic linear system, is asymptotically (and hence exponentially) stable if and only if there exists a periodic quadratic Lyapunov function verifying stability, i.e. 
\[
V(x_\perp,t) = x_\perp'P(t) x_\perp
\]
with $P(t)$ symmetric periodic positive definite (SPPD) such that 
\begin{equation}\label{eqn:lyap_ineq}
\dot P(t)+A(t)'P(t)+P(t)A(t)+H(t)\le 0
\end{equation}
for some SPPD $H(t)$. One example is the unique periodic solution of the Lyapunov differential equation:
\begin{equation}\label{eqn:lyap_eq}
\dot P(t)+A(t)'P(t)+P(t)A(t)+Q(t)=0
\end{equation}
for any $Q(t)=Q(t+T)>0$ which can be freely chosen. In other cases, if \eqref{eqn:sys} is the result of a linear control design procedure such as LQR or $H^\infty$, then $P(t)$ may come from solution of an associated Riccati equation.

\begin{theorem} \label{thm:lyapc} The following statements are equivalent
\begin{enumerate}
\item $x^\star$ is an exponentially orbitally stable $T$-periodic solution of \eqref{eqn:sys}. 
\item The transverse linearization \eqref{eqn:translin1} is exponentially stable.
\item For each SPPD $Q(t)$ there exists a unique SPPD solution $P(t)$ to the periodic Lyapunov equation \eqref{eqn:lyap_eq}.
Furthermore, for a sufficiently small region around $x^\star$ the Lyapunov function
\begin{equation}\label{eqn:lyap_func}
V(x) = x_\perp'P(\tau)x_\perp,
\end{equation}
verifies exponential orbital stability of $x^\star$.
\end{enumerate}
\end{theorem}

\begin{prf} 

\noindent $1 \rightarrow 2$: 
Consider dynamics in coordinates $(x_\perp, \tau)$. Due to exponential stability for a sufficiently small region around $x^\star$ dynamics are well posed and Lipschitz. The linearization of the dynamics about the target orbit has the form
\begin{equation}
\frac{d}{dt}\begin{bmatrix}\tau\\ x_\perp\end{bmatrix} \approx \begin{bmatrix}0 &\star \\ 0 & A(t)\end{bmatrix} \begin{bmatrix}\tau\\ x_\perp\end{bmatrix} 
\end{equation}
where $\star$ indicates a ``don't care'' element. This linearized system has a monodromy matrix (i.e. state transition matrix over the period $T$) of the form
\[
\Psi =  \begin{bmatrix}1 &\star \\ 0 & \Psi_\perp \end{bmatrix},
\]
where $\Psi_\perp$ is the monodromy matrix of the transverse linearization \eqref{eqn:translin1}. It is known that if $x^\star$ is exponentially orbitally stable, then $\Psi$ has one eigenvalue $\lambda_1=1$ and all other eigenvalues satisfy $|\lambda_i|<1$ (\cite{Hale80}). Due to the block diagonal form, this clearly implies that the eigenvalues of $\Psi_\perp$ are all inside the unit circle, which implies the stability of the transverse linearization.

\noindent $2 \leftrightarrow 3$: Is a known result in periodic linear systems theory (see, e.g.,  \cite{Bittanti91}). 

\noindent $3 \rightarrow 1$: We will use Lyapunov's direct method with the Lyapunov function \eqref{eqn:lyap_func} to verify exponential stability.

It is clear that $V(x)=0$ implies that $x_\perp=0$ hence $x\in \mathcal{X}^\star$. Furthermore, from \eqref{eqn:trans_full1}, \eqref{eqn:trans_full2}, and \eqref{eqn:lyap_ineq} we have
\begin{eqnarray}
\dot V(x)&\le& 2x_\perp P(\tau) h(x_\perp, \tau)+g(x_\perp, \tau)x_\perp'\frac{d}{d\tau}P(\tau)x_\perp \notag\\ &&\,-x_\perp H(\tau) x_\perp\notag
\end{eqnarray}
and $h=O(|x_\perp|^2)$ and $g=O(|x_\perp|)$, so all terms except $-x_\perp H(\tau) x_\perp$ are third-order in $x_\perp$, therefore for sufficiently small $\|x_\perp\|$ we have $\dot V(x)\le-\alpha\|x_\perp\|$ for some $\alpha>0$. 

Now, since $x^\star$ is a closed orbit, and by definition of $x_\perp$, over some compact region around $x^\star$ we have $k_1\dist(x,  \mathcal{X}^\star)\le \|x_\perp\| \le k_2\dist(x,  \mathcal{X}^\star)$ for some $k_1, k_2>0$, hence  $k_3\dist(x,  \mathcal{X}^\star)\le V(x) \le k_4 \dist(x,  \mathcal{X}^\star)$ for some $k_3, k_4>0$ and $\dot V(x) \le \alpha k_1\dist(x,  \mathcal{X}^\star)$ which proves exponential orbital stability.
\end{prf}

A simple bisection search can be done to find a region $V(x_\perp,\tau)\le \rho$ in which orbital stability is guaranteed. Then $V_1(x_\perp,\tau)= V(x_\perp,\tau)/\rho$ can be used as an initial seed in the search for more expressive Lyapunov functions around the orbit.

\section{Hybrid Systems}

In this section, we extend the above results to the case of hybrid systems, and propose an algorithm for solving Problem \ref{prob:regions_hybrid}.

Suppose we find a vector function $z(t)\in\mathbb R^n$ for $t\in[0, T)$ such that
\begin{eqnarray}
z(t)'f(x^\star(t)) &>&0, \forall \ t \in [0,T)\\
z(t_i)&=&\alpha_ic_i, \ \forall \ i =1, 2, ... N,
\end{eqnarray}
for some sequence of scalars $\alpha_i \ne 0$. That is, $z(t)$ aligns with the normals to the switching surfaces, and always makes a ``small angle'' with $f(x^\star)$. The non-grazing condition on impact surfaces ensures that many such vector functions exist.


With the coordinate system defined to line up with the switching surfaces, an entire region of the form $V(x_\perp,\tau)\le 1$ for a time $t_i$ will map into the image of $S_i$ under $\Delta_i$. 

The switching update of the transversal coordinates is:
\[
x_\perp^+=\Pi(\tau_i^+)\big[\Delta_i\big(x^\star(\tau_i^-)+\Pi(\tau_i^-)'x_\perp^-\big)-x^\star(\tau_i^+)\big]
\]
and the transverse linearization of the impact map is:
\[
A_d=\Pi(\tau_i^+)\frac{\partial\Delta_i}{\partial x}\Pi(\tau_i^-)'
\]
evaluated at $x=x^\star(\tau_i^-)$.

\subsection{Stability Conditions}

For systems with impacts, stability conditions will be of the form
\begin{eqnarray}
\frac{d}{dt}V(x_\perp, \tau)&<&0,\label{eqn:hybridc1}\\
z(\tau)'f(x^\star(\tau))-\frac{\partial z(\tau)}{\partial \tau}'\Pi(\tau)'x_\perp&>&0,\label{eqn:hybridc2}\\
\Big(c_{-}'\big(x^\star(\tau_i^-) +\Pi(\tau_i^-)'x_\perp\big)-d_{-}\Big)l_s(x_\perp, \tau) \notag&&\\+g(x^\star(\tau_i^-) +\Pi(\tau_i^-)'x_\perp\big)&<&0, \label{eqn:hybridc3}
\end{eqnarray}
for continuous phases over the regions  $V(x_\perp, \tau)\le 1$ via Lagrange multipliers. Here $l_s(x_\perp, \tau)$ is a Lagrange multiplier which is not constrained to be positive. The first constraint verifies stability via the decreasing Lyapunov function; the second verifies well-posedness of the change of variables; the third constraint ensures that the switching surface is not hit before it is expected. For some systems, it may be obvious that this will not happen, and the third condition could be dropped. At the switching times the conditions to be verified are:
\begin{eqnarray}
V\Big(\Pi(\tau_i^+)\big[\Delta_i\big(x^\star(\tau_i^-)+ \Pi(\tau_i^-)'x_\perp\big) &&\notag\\ -x^\star(\tau_i^+)\big], \tau_i^+\Big) -V(x_\perp, \tau_i^-)&\le&0, \label{eqn:hybridd1}\\
g\big(x^\star(\tau_i^-) +\Pi(\tau_i^-)'x_\perp\big)&>&0. \label{eqn:hybridd2}
\end{eqnarray}
The first verifies stability, the second that 
the state is within the region of definition of the impact on the switching hyperplane. All constraints evaluated over the regions $V(x_\perp, \tau)\le 1$ via Lagrange multipliers.

\begin{theorem}\label{thm:stabh}
Suppose conditions \eqref{eqn:hybridc1}--\eqref{eqn:hybridd2} hold in the region $D:=\{(x_\perp,\tau): V(x_\perp, \tau)\le 1\}$ then $D$ is an inner estimate of the region of attraction of the solution $x^\star(\cdot)$ of the hybrid system \eqref{eqn:sysh1}, \eqref{eqn:sysh2}.
\end{theorem}
The proof is omitted because of space restrictions but is a minor variation on the proof of Theorem \ref{thm:stab} and standard application of Lyapunov's direct method to impulsive systems (see, e.g., \cite{Bainov89}).

\subsubsection{Sum-of-Squares Relaxation}

The above conditions for orbital stability of hybrid limit cycles can be relaxed to sum-of-squares programs analogously to \eqref{eqn:sos1}--\eqref{eqn:sos4} with additional SoS constraints added for the impact conditions \eqref{eqn:hybridd1}, \eqref{eqn:hybridd2}. Details are omitted due to space restrictions.

\subsection{Initial Candidate Lyapunov Function}

A natural candidate Lyapunov function for the hybrid case is the unique SPPD solution of the jump-Lyapunov function:
\begin{eqnarray}
-\dot P &=& A'P+PA +Q, \ t \ne t_i \notag\\
P(\tau_i^-)&=&A_d(\tau_i)'P(\tau_i^+)A_d(\tau_i)+Q_i, \ t = t_i \notag
\end{eqnarray}
A similar statement to Theorem \ref{thm:lyapc} can be made applying the results of \cite{Nicolao98} and \cite{Bainov89}. Details are omitted due to space restrictions.

\section{Controlled Systems and Orbital Stabilization}

Consider a controlled system:
\begin{equation}\label{eqn:sys_u}
\dot x(t) = f(x(t),u(t))
\end{equation}
with $x(t)\in\mathbb R^n$ and $u(t) \in \mathbb R^m$.
Suppose a nominal $T$-periodic solution has been found: $x^\star(t), u^\star(t)$.  We wish to study perturbations and orbital stabilization of this system. With the nominal control implemented as a function of time, the resulting perturbed system:
\[
\dot x(t) = f(x^\star(t)+\Delta x(t), u(t))
\]
is non-autonomous has extremely complicated dynamics and orbital stability cannot be characterized locally.

In contrast, if $u(t)$ is a feedback policy on $x(t)$ then the region-of-attraction analysis reduces to the study of an autonomous system and the methods above can be applied directly. This analysis can of course be applied with any feedback control law, however the analysis of the transverse dynamics suggests a natural approach to the stabilization problem.

Suppose a set of transversal surfaces is chosen, along with a projection $\bar\tau(\cdot): \mathbb R^n \rightarrow [0, T)$ such that $x\in S(\bar\tau(x))$ for all $x$ in a neighbourhood\footnote{This neighbourhood will be characterized by the computed regions of attraction of the closed-loop system, but is guaranteed to exist locally which is sufficient for the present arguments.} of the solution $x^\star(\cdot)$, and we apply the control 
\[
u(x) = u^\star(\bar\tau(x))+\bar u.
\]
with $\bar u$ a stabilizing control, for example,
\[
\bar u = -K(\bar\tau(x))\Pi(\bar\tau(x))[x-x^\star(\bar\tau(x))]
\]
where $K(\cdot):[0,T)\rightarrow \mathbb R^{m\times n}$ is feedback gain. Note that the nominal command $u^\star$ is being applied as  feedback law based on the estimated phase $\tau$ rather than as an open-loop function of time.

This feedback gain $K(\cdot)$ can be designed using the controlled transverse linearization:
\begin{equation}\label{eqn:translin_u}
\dot x_\perp = A(t)x_\perp(t) +B(t)\bar u(t), \ \ t\in [0,T)
\end{equation}
with $A(t)$ is in \eqref{eqn:translinA} with $f(x^\star(\tau)$ replaced by $f(x^\star(\tau),u^\star(\tau))$ and
 $\frac{\partial f(x^\star(\tau))}{\partial x}$ replaced by $\frac{\partial f(x^\star(\tau),u^\star(\tau))}{\partial x}$. The formula for $B(t)$ is
\begin{equation}\label{eqn:translinB}
B(\tau) = \Pi(\tau)\frac{\partial f(x^\star(\tau),u^\star(\tau))}{\partial u}-\Pi(\tau)f(x^\star(\tau), u^\star(\tau))\frac{\partial \dot \tau}{\partial u}
\end{equation}
where
\[
\frac{\partial \dot \tau}{\partial u}=\frac{1}{z(\tau)'f(x^\star(\tau), u^\star(\tau))}z(\tau)' \frac{\partial f(x^\star(\tau),u^\star(\tau))}{\partial  u}.
\]
Note that if orthogonal transversal surfaces are used then $\Pi(\tau)f(x^\star(\tau),u^\star(\tau))$ is zero and therefore the second term in \eqref{eqn:translinB} is zero.

Any $K(\cdot)$ which stabilizes the periodic system \eqref{eqn:translin_u} locally exponentially stabilizes the target orbit. An associated Lyapunov function for the linear system can be used for regional verification of the closed-loop orbital stability as above.

\subsection{LQR via Jump-Riccati Equation}

As an example of possible control design procedure, we give the solution to the transverse-LQR problem via a jump Riccati equation. This solution is far from the only one, but is quite generic in that if the orbit is stabilizable in a certain sense, this method will always stabilize it.

\begin{theorem}
The target orbit $x^\star(\cdot)$ of the original nonlinear system is exponentially orbitally stabilizable by smooth state feedback if and only if the associated transverse linearization is stabilizable.
Furthermore, if this is the case then there exists a unique SPPD solution of the following jump-Riccati equation:
\begin{eqnarray}
-\dot P &=& A'P+SA -PBR^{-1}B'P +Q, \ t \ne t_i \notag\\
P(\tau_i^-)&=&A_d(\tau_i)'P(\tau_i^+)A_d(\tau_i)+Q_i, \ t = t_i \notag
\end{eqnarray}
with $R, Q, Q_i > 0$. An associated locally stabilizing feedback controller for the original nonlinear system is then given by 
\[
u(\tau, x_\perp) = u^\star(\tau)-R^{-1}B(\tau)'P(\tau)x_\perp.
\]
\end{theorem}

\begin{prf}
The full proof is omitted due to space restrictions, however it follows the equivalence to exponential stabilization via smooth state feedback and stabilization of the transverse linearization from straightforward application of Theorem \ref{thm:stabh}. The existence of the SPPD solution of the jump-Riccati equation follows from minor modification of the stabilizability properties of periodic impulsive linear systems \cite{Nicolao98, Bittanti91}. Note that the restriction that $Q, Q_i>0$ can be relaxed somewhat to a discrete-continuous observability property of the cost function.
\end{prf}


\section{Optimization of Transversal Surfaces}\label{sec:choose_z}

The method described above allows great flexibility in the choice of transversal surfaces, parameterized by their normal vectors $z(\tau)$. In some cases, problem specific information might suggest a natural candidate (see examples). In others, the classical choice $z(\tau)=f(x^\star(\tau))/\|f(x^\star(\tau))\|$ may be sufficient. However, it is likely that for many practical examples some optimization of $z(\tau)$ is appropriate.

The well-posedness conditions break when
\[
z(\tau)'f(x^\star(\tau))-\frac{\partial z(\tau)}{\partial \tau}'\Pi(\tau)'x_\perp=0.
\]
It is natural to try to choose $z(\tau)$ such that this is prevented from happening in a large region around the orbit, i.e. for a large ball in $\|x_\perp\|$.

Finding the smallest $\|x_\perp\|$ for which this can happen is simply a linear-constrained least-squares problem, solved with the augmented cost function
\[
J(x_\perp, \lambda) = \frac{1}{2}x_\perp'x_\perp+\lambda\Big(z(\tau)'f(x^\star(\tau))-\frac{\partial z(\tau)}{\partial \tau}'\Pi(\tau)'x_\perp\Big)
\]
The usual computation leads to the following minimizer:
\[
x_\perp^\dagger=\frac{z(\tau)'f(x^\star(\tau))\Pi(\tau) \frac{\partial z(\tau)}{\partial \tau}}{\left|\Pi(\tau) \frac{\partial z(\tau)}{\partial \tau}\right|^2}
\]
and since we only care about the size of $x_\perp$, this simplifies to
\[
|x_\perp^\dagger|=\frac{|z(\tau)'f(x^\star(\tau))|}{\left|\Pi(\tau) \frac{\partial z(\tau)}{\partial \tau}\right|}.
\]
Note that $\frac{\partial z(\tau)}{\partial \tau}$ is orthogonal to $z(\tau)$ and hence entirely in the row-space of $\Pi(\tau)$, so $|\Pi(\tau) \frac{\partial z(\tau)}{\partial \tau}|=|\frac{\partial z(\tau)}{\partial \tau}|$ and
\[
|x_\perp^\dagger(\tau)|=\frac{|z(\tau)'f(x^\star(\tau))|}{|\frac{\partial z(\tau)}{\partial \tau}|}
\]
making it unnecessary to compute $\Pi(\tau)$ at every step.


The size of $x_\perp^\dagger(\tau)$  should be maximized in some sense. However, notice that the denominator has $\frac{\partial z(\tau)}{\partial \tau}$ as a factor. The optimum may have $z(\tau)$ constant for some intervals of $\tau$ which means $|x_\perp^\dagger(\tau)|$ would go to infinity, which may make an optimization difficult. It is therefore better-posed to minimize its inverse. E.g. one can minimize
\begin{eqnarray}
z(\tau)^\star &=& \arg\min\left(\int_0^T \frac{|\frac{\partial z(\tau)}{\partial \tau}|^p}{|z(\tau)'f(x^\star(\tau))|^p} d\tau\right)^{1/p}\\
&\textrm{s.t.}& \ z(\tau)'f(x^\star(\tau))>\delta\\
&& \ \|z(\tau)| = 1 \ \forall \, \tau \in [0,T].
\end{eqnarray}
for some small positive $\delta$. 
The author has had success choosing $p$ to be quite large, say between 10 and 100, so that it approximates an $L^\infty$ norm but remains a smooth optimization problem.

For systems with impacts, the integrand in the above optimization can be multiplied by a smooth shaping function $\phi(\tau)$ with $\phi(\tau_i)=0$, since the transversal surfaces must be aligned with the switching surfaces and cannot be optimized.

This is a nonconvex optimization, so it requires a good initial seed. A natural candidate optimization is $z(\tau)=f(x^\star(\tau))/\|f(x^\star(\tau))\|$. This initial guess always satisfies the transversality conditions, but may not always satisfy the condition of alignment to switching surfaces.

\section{Illustrative Examples}

The proposed method has been used by the author and colleagues to compute regions of stability for a number of example systems. The systems so far considered include: the van der Pol oscillator, which illustrates the importance of selection of $z(\tau)$; the rimless wheel, a very simple model of an underactuated walking robot that illustrates the extension to hybrid systems; the compass-gait walker, a more complex model of underactuated walking with regions of stability that cannot be found in closed form, and illustrate the proposed method's use for nontrivial systems. A detailed discussion of these examples can be found in the companion paper \cite{Manchester10b}.

\section{Conclusions and Future Work}

We have presented a constructive procedure for computation of an inner (conservative) estimate of the region of attraction of a nonlinear hybrid limit cycle. Such a procedure was motivated by problems of control and motion planning of walking robots, but undoubtedly will find applications in many other fields where oscillating nonlinear systems are of interest.

In the case where a feasible but unstable periodic solution is known, the method of transverse coordinates can also be used to construct a feedback controller which guarantees exponential orbital stability. The regions of attraction of the closed-loop system can then be analyzed using the proposed technique.

The method of analysis via Lyapunov functions and sum-of-squares relaxation has a natural extension to uncertain systems via storage functions and integral quadratic constraints (see, e.g., \cite{Megretski97, Petersen00}) or parameter-dependent Lyapunov function (\cite{Topcu10}). Robustness of the existence of periodic trajectories has previously been studied using an extension of the small-gain theorem in 
\cite{Jonsson05, Jonsson08}. Future work will include connections between such robust analysis tools and the present work on regions of stability.

\bibliographystyle{IEEEtran}
\bibliography{elib}

\end{document}